\documentclass{amsart}

\newtheorem{theorem}{Theorem}[section]
\newtheorem{lemma}[theorem]{Lemma}
\newtheorem{proposition}[theorem]{Proposition}
\theoremstyle{definition}
\newtheorem{definition}{Definition}[section]
\theoremstyle{remark}
\newtheorem{remark}{Remark}[section]
\newtheorem{example}{Example}[section]
\numberwithin{equation}{section}
\newtheorem{corollary}{Corollary}[section]

\begin{document}

 \title[Filling growth]{The ends of manifolds with bounded
   geometry, linear growth and finite filling area}

\author[L.Funar]{Louis Funar}
\address{Institut Fourier BP 74, UMR 5582 CNRS, Universit\'e de  
Grenoble I, 38402 Saint-Martin-d'H\`eres cedex, France}
\email{funar@fourier.ujf-grenoble.fr}  
\thanks{Partially supported by GNSAGA and MIUR of Italy.} 
\thanks{Preprint available at {\tt 
http://www-fourier.ujf-grenoble.fr/\~\,funar}}

\author[R.Grimaldi]{Renata Grimaldi}
\address{Dipartimento di Matematica ed Applicazioni, Facolta di
  Ingegneria, Universit\`a di Palermo, Viale delle Scienze, 
90128, Palermo, Italia}
\email{grimaldi@math.unipa.it}

\date{}
 
 \begin{abstract} 
We prove that simply connected open manifolds 
of bounded geometry, linear growth and sub-linear filling growth 
(e.g. finite filling area) are simply connected at infinity. 

\noindent {\em MSC}:  {53 C 23, 57 N 15.}



\end{abstract}



\maketitle
 
 \section{Introduction}
A ubiquitous theme in Riemannian geometry is the relationship between
the geometry (e.g. curvature, injectivity radius) and the topology. In
studying non-compact manifolds constraints come from the asymptotic
behaviour of geometric invariants (e.g. curvature decay, volume
growth) as functions on the distance from a base point.  The expected 
result is the manifold tameness out of geometric
constraints. This is illustrated by the classical theorem of
Gromov which asserts that a complete hyperbolic manifold of finite volume and 
dimension at least 4 is the interior of a compact manifold with
boundary.  Our main result below yields tameness in the case  
when the filling area is finite, for those manifolds having bounded geometry
and linear growth. We recall that: 

\begin{definition}
A non-compact Riemannian manifold has  
{\em bounded geometry} if the injectivity radius $i$ is bounded from
below and the absolute value of the curvature $K$ is bounded from above.  
\end{definition}
\begin{remark}
One can rescale the metric in order that $i\geq 1$ and $\mid K\mid
\leq 1$ hold. 
\end{remark}

\begin{definition}
The {\em filling  area function} $F_X(l)$ of the simply connected  
manifold $X$ is the smallest number with the property that any loop of length 
$l$ bounds a disk of area $F_X(l)$. 
\end{definition}

\noindent It is customary to introduce the following equivalence relation: 
\begin{definition}
Two positive real functions are {\em equivalent}, and one writes
$f\sim g$, if 
\[ c_1f(c_2 x) + c_3 \leq g(x) \leq C_1 f(C_2 x) + C_3, \]
for positive $C_i, c_i$. By abuse of language we will call 
{\em filling area} the equivalence class of the filling area function. 
\end{definition}

\begin{remark} 
The equivalence class of the  filling area function $F_{\widetilde{X}}$
of the universal covering space of a  
compact manifold  $X$ is independent of the metric we chose on 
$X$. In fact the filling area is a quasi-isometry invariant, 
and hence an invariant  of the fundamental group of $X$, which we
keep calling the {\em filling area of the group}. 
\end{remark}

\noindent There are two interesting classes of open Riemannian manifolds of
bounded geometry:
\begin{itemize} 
\item the universal coverings of compact manifolds. Up to
  quasi-isometry these manifolds are determined by their deck
  transformations groups. This is part of Gromov's program of
  classifying discrete groups up to quasi-isometry. 
\item manifolds of bounded geometry with linear volume growth. 
It is known that for any 
  $\varepsilon > 0$ there exists a metric on a given open manifold 
which has bounded geometry and growth less that $1+\varepsilon$. 
Thus super-linear growth seems to be topologically unobstructed. 
However there exists (see \cite{FG,G}) a complete 
topological characterization of those 
manifolds supporting a metric of linear growth and bounded geometry, 
following Cheeger and Gromov (\cite{CG}): these are the 
manifolds of finite topology at infinity. Specifically, this means
that there exists a proper Morse function $\lambda$, 
such that the level hypersurfaces $\lambda^{-1}(n)$ 
are  pairwise diffeomorphic for all $n\in {\bf Z}$.    
\end{itemize}

\noindent The filling area has been studied for universal coverings and
shown to be equivalent to the Dehn function of the group, measuring
the complexity of the group.
This topic received recently a lot of consideration (see \cite{OS} for a
survey). 
\begin{remark}
It is worth mentioning that sub-quadratic 
filling area implies linear filling area and this is equivalent to
the group being word-hyperbolic. On the other hand the exponents of
the  polynomial filling areas fill in a dense subset in $[2,\infty)$.   
\end{remark}
 
\noindent This motivates the study of the filling function in the 
second case, as well. 
The definition of the filling area can be quite
unappropiate for general non-compact Riemannian manifolds. 
Even for manifolds of bounded geometry, it  
might take only infinite values, in  which case it does not give 
any valuable information about the topology. 
We introduce for this reason the  following refined version: 

\begin{definition}
Let $X$ be non-compact. 
The {\em filling area function} $F_X(l, r)$ is the smallest area
of the disk in $X$ filling a loop of length $l$ lying 
in the metric ball $B_X(r)$ of radius $r$ on $X$. 
\end{definition}

\begin{definition}
The positive functions $f(l,r)$ and $g(l,r)$ are equivalent, and  we 
write $f\sim g$, if 
\[ c_1(l)f(c_2(l), r) + c_3(l) \leq g(l,r) \leq C_1(l) f(C_2(l),r) + C_3(l), \]
for positive increasing functions $C_i(l), c_i(l)$. By 
{\em filling (area) growth} one means the equivalence class of the
filling  area function. 
\end{definition}

\begin{remark}
The filling growth is a quasi-isometry invariant. 
\end{remark}

\begin{definition}
The filling growth  is said to be {\em sub-linear} if: 
\[\lim_{r\to \infty} \frac{F_X(l, r)}{r} = 0, \]
holds for every $l$. 
\end{definition}

\begin{remark}
If the filling area is finite then it is automatically sub-linear. 
\end{remark}

\begin{definition}
A non-compact polyhedron  $X$ is {\em simply connected at infinity} 
(s.c.i.), and we write also  $\pi_1^{\infty}(X)=0$, if given   
a compact set $K\subset X$ there exists another compact set $L$ with 
$K\subset L\subset X$, such that any loop in $X-L$ is null-homotopic in 
$X-K$. Alternatively the map induced by the inclusion 
$\pi_1(X-L)\to \pi_1(X-K)$ is trivial (i.e. the zero map).  
\end{definition}
\begin{remark}
Some authors call this $\pi_1$-{\it triviality at 
infinity} or {\it 1-LC at infinity} and reserve the term 
s.c.i. for the special case in which $L$ can be chosen so that, in
addition, $X-L$ is connected. These notions are equivalent for one
ended spaces, such as contractible spaces. 
\end{remark}

\begin{remark}
The simple connectivity 
at infinity is an
important {\em tameness} condition on the ends of the space.
It has been used to characterize Euclidean space 
among contractible open topological 
$n$-manifolds by Siebenmann (\cite{Sie})  for $n\geq
5$, Freedman (\cite{Free}) for $n=4$ and by Edwards
(\cite{Ed}) and Wall (\cite{Wall})
 for $n=3$ (after assuming the irreducibility to avoid the 
Poincar\'e conjecture). An earlier related result is the Stallings-Zeeman 
engulfing theorem (see \cite{Sta}),  one of whose consequences is the fact 
that an open contractible PL $n$-manifold ($n\geq 5$) which is s.c.i. 
(notice that Stallings used  a slightly stronger notion of s.c.i. than 
that commonly used now) 
is PL-homeomorphic to the Euclidean
space.  It is thus of some
interest of finding criteria which imply that a space is s.c.i.
\end{remark}

\noindent The main result of our paper is the following topological 
characterization: 

\begin{theorem}
A simply connected open Riemannian manifold of bounded geometry, linear growth
and sub-linear filling growth is simply connected at
infinity. 
\end{theorem}

\begin{corollary}
An open contractible Riemannian manifold of bounded geometry, linear growth
and sub-linear filling growth is diffeomorphic 
(only homeomorphic in dimension 4) to the Euclidean space 
(one needs to assume the irreducibility in dimension 3).  
\end{corollary}

\begin{remark}
The s.c.i. for universal coverings is a group theoretical property, 
which is not yet fully understood. 
A.~Casson, V.~Poenaru and M.~Mihalik, among others, 
gave geometric conditions on the group insuring 
the s.c.i. Nevertheless these are not
general enough in order to include all 
infinite 3-manifold groups, as claimed by the uniformization conjecture. 
\end{remark}

\noindent {\bf Acknowledgements.} The authors are indebted to Pierre
Pansu  and the referee for helpful comments and advice.

\section{A geometric finiteness result}
Before we proceed let us fix the notations. We denote by $K_M$, $i_M$ 
the sectional curvature and the injectivity radius of the 
Riemannian manifold
$M$. When $W\subset M$ is a submanifold, $i_{W\subset M}$ stands for
the normal injectivity radius (i.e. the maximal radius of a tube 
around $W$  which is  embedded in $M$).  By ${\rm II}_{W}$ 
one denotes the second fundamental 
form of $W$. We fix once for all a base point $p$ on $M$, so that all 
metric balls in the sequel   
are centered at $p$, unless the opposite is explicitly stated. 
The Riemannian metric on $M$ induces a distance 
$d_M$ (or $d$ if no confusion arises), a volume form $vol$ and 
a 2-dimensional $area$. 
The length of the curve $\gamma$ is usually denoted by
$l(\gamma)$.

The result of this section adds some  geometric control
to the finiteness theorem of \cite{FG}. Specifically, we will prove below
that:   
\begin{proposition}\label{lip}
Let $M$ be an open simply connected Riemannian manifold of bounded geometry 
and linear growth. There exist then:
\begin{enumerate}
\item  an exhaustion of $M$ by compact
submanifolds $M_j\subset M_{j+1}\subset ...\subset M$,  
\item tubular neighborhoods $Z_j$ of $\partial M_j$ in $M$, 
\item constants $c$, $A$ and $\delta >0$ depending only on $M$,  
\item a closed Riemannian manifold $V$ of dimension one less than $M$, 
\end{enumerate}
such that the Lipschitz distance between 
$Z_j$ and $V\times [0,c\delta]$ is bounded by $c$, and also
$d(p, \partial M_j)\leq Aj$, for all $j$.  
\end{proposition} 
\begin{proof}
We need first a slight improvement of a lemma due to Cheeger and
Gromov: 
\begin{lemma}\label{chop}
Let $M$ be a Riemannian manifold of bounded geometry, 
$X\subset M$ a subset, $\varepsilon > 0$. We denote by 
$T_{\varepsilon}(X)=\{ x; d(x, X) < \varepsilon\}$ the set of points
having distance less than $\varepsilon$ from $X$.  
There exists a hypersurface $W\subset T_{\varepsilon}(X)-X$ such that:
\begin{enumerate}
\item  $vol(W)\leq \alpha(\varepsilon)(vol(T_{\varepsilon}(X)-X))$. 
\item  $\parallel {\rm II}_{W}\parallel\leq \beta(\varepsilon)$.
\item $i_W\geq c(\varepsilon)> 0$. 
\item $i_{W\subset M}\geq c(\varepsilon)> 0$. 
\end{enumerate}
where $\alpha, \beta, c$ denote constants depending only on $\varepsilon$. 
\end{lemma}
\begin{proof}
The first two requirements are granted by the approximation
theorems from (\cite{CG1}, p.127-135) and (\cite{CG}). 
The third part was proved in lemma 3 from \cite{FG}. 
It remains therefore to prove that the hypersurface $W$ provided there can be 
supposed to satisfy also the fourth condition.

Moreover, the normal injectivity radius can be  controlled by means of
the second fundamental form and the curvature: 
\begin{lemma}\label{doi}
If $\parallel {\rm II}_{W}\parallel\leq A$ and $K_M\leq A$, then
$i_{W\subset M}\geq c(A)> 0$ holds true.  
\end{lemma}
\begin{proof}
This result seems to be well-known. One has to 
bound from below the distance to the closest focal point.
The Jacobi theorems hold for the index form associated to the compact
submanifold $W$ (see e.g. \cite{Ch,CE}): thus the normal unit 
speed geodesic $\gamma$
has no conjugate points in some interval iff the index form is 
positive. The proof of the Morse-Schoenberg (or Rauch comparison)
theorem extends (see \cite{Ch}, p.79) 
without essential modifications to this setting 
(see also \cite{P}, p.172-175). 
Otherwise one can use the more general Rauch comparison theorem
sketched in (\cite{Ch}, Remark 3.4., p.135). 
Effective estimates for the constant are given in \cite{H} 
where only positively curved manifolds are considered, 
but the arguments work for curvature bounded from above. 
\end{proof}
\noindent This yields the claim of lemma \ref{chop}.
\end{proof}

\noindent A result similar to lemma \ref{chop} has been 
used in \cite{FG}, in order to 
obtain the existence of an exhaustion $M_j$ whose 
boundaries $\partial M_j$ belong to a  finite family of 
diffeomorphism types. One has then to see that the Lipschitz
distances are uniformly bounded and further that the distance 
to the boundaries grows linearly.

The lemma \ref{chop} yields uniform lower bounds (away from zero)
for the normal injectivity radii $i_{\partial M_j\subset M}$, 
and  the fact that the  metrics induced on 
the manifolds $\partial M_j$ (which are pairwise diffeomorphic) 
are at uniformly bounded Lipschitz distance one from each other. 
This implies the first assertion from proposition \ref{lip}.

In what concerns the second claim of proposition \ref{lip}, we  follow 
closely the proof  of the main result from \cite{FG}. 
Specifically we can state:  
\begin{lemma}\label{cost}
There exists an exhaustion $\{ M_j\}$ and constants $c_i=c_i(M)$ such that:  
\begin{enumerate}
\item the boundary manifolds 
$\partial M_j$ inherit metrics fulfilling 
\[ \mid K_{\partial M_j}\mid \leq c_1, \; 
i_{\partial M_j}\geq c_2 > 0, \; vol(\partial M_j)\leq c_3, \mbox{ for
  all } j,\]
\item the boundary $\partial M_j$ is sandwiched between two metric
  spheres, namely 
\[ c_4 j\leq d(p, \partial M_j)\leq c_5 j,  \mbox{ for
  all } j.\]
\end{enumerate}
\end{lemma}
\begin{proof}
In order to obtain (2) it suffices to improve the 
sub-lemma 2.1 from (\cite{FG}, p.853), as follows: 
\begin{lemma}
Let $A(n)=B_M(n+1)-B_M(n)$ be the annuli of unit width, where 
$B_M(n)$ denotes the ball of radius $n$ centered at $p$. 
Then there exists a sequence  ${n_k}$ and constants $c, C$ such that 
\begin{enumerate}
\item $vol(A(n_k))\leq c$, 
\item $n_k\leq C k$. 
\end{enumerate}
\end{lemma}
\begin{proof}
The linear growth hypothesis is that  $vol(B_M(n))\leq a n$, 
for some constant $a$. Let $S=\{n\in \mathbb{Z}_+; vol(A(n))\leq 2a\}$. 
Assume that $card(S\cap \{ 1,2,...,N\}) <\frac{N}{2}$. 
Then $vol(B_M(N)) > (2a)\frac{N}{2}=aN$, contradiction.
This shows that one can take $c=C=2$ above. 
\end{proof}  
\noindent In particular there exists an exhaustion $M_j$  with 
the property that all $\partial M_j$ are diffeomorphic 
and moreover $d(p, \partial M_j)\leq AC j$, where $A$ is  
the number of diffeomorphism types of manifolds  
satisfying the first condition in lemma \ref{cost} (and $A$ being
finite by Cheeger's finiteness theorem). This proves lemma \ref{cost}.
\end{proof}
\noindent As already remarked this lemma ends the proof 
of proposition \ref{lip}. 
\end{proof}

\section{The proof of the theorem}
Assume that we fixed an exhaustion $M_j$ and the tubular neighborhoods 
$Z_j$  like in the proposition \ref{lip}. One knows that $Z_j$ is canonically
identified with a cylinder $\partial M_j\times [0, r_j]$, 
where all $r_j$ are bounded from below by some $\delta >0$. 

\begin{lemma}\label{disk}
There exists a constant $c_6$ with the following property. 
For any  2-disk $D$  transverse to $Z_j$, such that  
$D\cap Z_j$ is not null-homotopic (i.e. 
$D\cap M_j\times\{ t\}$ is not null-homotopic, for any $t\in[0,r_j]$)
we have  $area(D\cap Z_j)\geq c_6$. 
\end{lemma}
\begin{proof}
Let $cl_0$ be the 1-systole of $V$, i.e. the 
length of its smallest closed geodesic. 
By proposition \ref{lip} the length of each 
non-trivial component of $D\cap \partial M_j\times \{t\}$ 
is bounded from below by $l_0$.
Then  the projection $p:\partial M_j\times[0,\delta]\to [0,\delta]$ is  
decreasing the length, i.e. $\parallel\nabla p\parallel \leq 1$. 
We take then $c_6=\delta l_0$. In fact the coarea formula states that:
\[ area(D\cap \partial M_j\times[0,\delta])=
\int_0^ \delta d t \int_{D\cap \partial M_j\times \{t\}} 
\frac{d s}{\parallel \nabla p\parallel} \geq \]
\[\geq \int_0^\delta d t \int_{D\cap \partial M_j\times \{t\}} d s
\geq 
\int_0^\delta l({D\cap \partial M_j\times \{t\}})d t \geq \delta l_0 
\]
\end{proof}

\noindent {\em Proof of the theorem.} 
Suppose that the contrary holds so that 
$M$ is not simply connected at infinity. 
\begin{lemma}
There exists
then a  compact $K$ such that the maps induced by the inclusions
$\iota_j:\pi_1(\partial M_j)\to \pi_1(M-K)$ are non-zero, for all
large enough $j$. 
\end{lemma}
\begin{proof}
There exists a compact $K$ such that arbitrary far loops bound only
disks touching $K$. If $K\subset M_{j_0}$ consider loops 
$l_j$ outside $M_j$ with this property. Any disk $D_j$ bounding $l_j$ should 
intersect $K$ hence $\partial M_j$. One puts the disk $D_j$ in general 
position with respect to $\partial M_j$. It follows that at least 
one loop component of $\partial M_j\cap D_j$ is not null-homotopic. 
Further the image of at least one loop component  
in $\pi_1(M-K)$ is non-zero. In fact, otherwise $l_j$ would bound 
a disk not touching $K$, when capping off the null-homotopy disks 
of these loops  with the annuli $D_j\cap M-M_j$. This contradicts the
choice of $l_j$. 
\end{proof}  

\noindent Recall that  all  
$\partial M_j$ are diffeomorphic to the closed manifold $V$. 
If  $a_1,...,a_N$ is a system of generators of $\pi_1(V)$, one has 
then an induced system  of generators  for $\pi_1(\partial M_j)$, which we
keep denoting by the same letters. There exists therefore
$k\in\{1,2,...,N\}$, such that the images $\iota_j(a_k)$ 
are non-zero, for all large enough  $j$.  
Let now $a_{k,j}$ denote a loop in $\partial M_j$ 
representing $a_k$ in $\pi_1(\partial M_j)$. 

\begin{lemma}
There exists some $c_7$ such that representatives loops $a_{k,j}$ 
can be chosen to be of length uniformly bounded by $c_7$, for all 
$k$ and $j$. 
\end{lemma}
\begin{proof}
This is a consequence of the fact that $\partial M_j$ are at uniformly bounded
Lipschitz distance from $V$. 
\end{proof}

\noindent Consider now a disk $D_j$ in $M$, filling the loop $a_{k,j}$. 

\begin{lemma}\label{ar}
For large enough $j$ we have $area(D_j)\geq \frac{c_6}{2}j$. 
\end{lemma}
\begin{proof}
One can assume that $D_j$ is transverse to all $Z_i$'s. 
By hypothesis $D_j\cap K\not=\emptyset$. Let $j_0$ be the 
minimal number with the property that
$K\subset M_{j_0-1}$. We claim that for all $i\in\{j_0,...,j\}$ 
the intersections $Z_i\cap D_j$ are not null-homotopic. 
Otherwise we could replace the disk $D_j$ by a more economical one, 
contained in $M_j - M_{j_0}$, which would contradict our choice of $j_0$. 
Further, by using lemma \ref{disk} one derives: 
\[ area(D_j) \geq \sum_{i=j_0}^{j} area(D_j\cap Z_i) \geq c_6(j-j_0).\]
This proves the lemma \ref{ar}.\end{proof}
 
\noindent Eventually we observe that the loops $a_{k,j}$ have 
uniformly bounded length and stay at distance  $c_5j$ far from the 
base point, because $a_{k,j}\subset \partial M_j$. 
In particular the filling growth  is linear in the
radius. This contradicts our hypothesis, and hence the theorem follows.

\section{A counterexample to the converse}
The metrics of bounded geometry on open s.c.i. manifolds
need  not  have a sub-linear filling area, as the following 
example shows: 
\begin{example}
There exists a Riemannian  metric $g$ of bounded geometry on 
${\bf R}^n$ (for $n\geq 3$) such that 
\[\lim_{r\to \infty} \frac{F_{({\bf R}^n,g)}(l, r)}{r} \neq 0.\]
\end{example} 
\begin{proof}
The method consists in modifying a decomposition of ${\bf R}^n$
into compression bodies by adding trivial cylinders of sufficiently 
large length. A loop bounding in a compression body can be
translated along the cylinder and thus the  bounding disk for the new loop
is the union of the former with a long cylinder. In particular we can
achieve a linear area function for this type of loops.   

Specifically let us focus on $n=3$, and set $T_j$ for an increasing 
family of solid tori, such that:
\begin{enumerate}
\item $T_{j+1}-int(T_{j})$ is a product $T^2\times [0,1]$ for all $j$ but an 
infinite sequence $A=\{j_1,j_2,...,j_k,...\}\subset \mathbb Z_+$. 
Here $T^2$ stands for the 2-torus. 
\item if $j\in A$ then the inclusion $T_j\hookrightarrow T_{j+1}$ is 
trivial up to isotopy, and thus $T_j$  is contained in a ball
$B^3$ embeded in $T_{j+1}$. Hence   $T_{j_i+1}-int(T_{j_i})=V^3$ has 
a fixed diffeomorphism type, and it is not diffeomeorphic to
$T^2\times [0,1]$. In particular the map
$\pi_1(\partial T_{j_{i}+1})\to \pi_1(T_{j_{i}+1}-int(T_{j_i}))$
is  not zero. 
\end{enumerate}
We consider a metric structure $g$ on each $T_{j+1}-int(T_{j})$, which is
a product  along the boundary (corresponding to each one of the two models, 
$V^3$ or $T^2\times [0,1]$), with isometric boundaries at unit
distance from each other.
We obtain a
metric on ${\bf R}^3-T_1$, which can be completed by capping off with
a Riemannian structure on $T_1$. 
This metric has linear growth. 

Let $c$ be a loop on $\partial T_{j_{i}}$, which is not
null-homotopic and does not bound in 
${\bf R}^3-T_{j_{i}}$, and thus it is not a longitude. 
Then any disk $D^2$ filling $c$ has a component in $T_{j_{i}}$, and
the innermost circle component of $D^2\cap \partial T_{j_{i}}$ is 
still homotopically non-trivial. We can suppose that 
this component is actually  $c$, and hence
that $c$ bounds in $T_{j_{i}}$. It follows that
$D^2\not\subset T_{j_{i}}-T_ {j_{i-1}+1}$, because 
$T_{j_{i}}-T_ {j_i+1}$ is a cylinder $T^2\times [0, j_i-j_{i-1}]$ 
retracting on $\partial T_{j_i}$. Therefore $D^2$ intersects 
$T_{j_{i-1}}$ and thus $D^2$ intersects non-trivially (i.e. along 
loops which are not null-homotopic) the $j_{i}-j_{i-1}$ 
intermediary tori $\partial T_j$ with $j_{i-1}+1\leq j\leq j_i$. Since these 
tori have uniformly bounded normal injectivity radii (thus  neighborhoods
isometric with $T^2\times [0,\delta]$) we obtain that 
$area(D^2)\geq \delta (j_{i}-j_{i-1})$. Now $\partial T_j$ is  
at distance $j-1$ from $T_1$ and thus we find  that: 
\[ \lim_{r\to \infty}\frac{F_{({\bf R}^3,g)}(l, r)}{\delta r}
\geq \lim_{i\to \infty}\frac{j_{i}-j_{i-1}}{j_i}. \]
By choosing $j_i$ growing fast enough we can insure that the right hand 
limit is $1$. 
\end{proof}

\section{Comments}
\noindent There exists a further refinement of the 
filling area function, 
as follows. Set $f_X(l,r;\lambda)$ for the smallest area of the disk 
in $X$ filling an arbitrary loop of length $l$ lying in the annulus
$B_X(r)-B_X(\lambda r)$ of $X$.  The function
$f_X(l,r;\lambda)$ need not being increasing anymore. 
The filling  growth is said to be weakly sub-linear if 
\[\lim_{r\to\infty}\inf \frac{f_X(l,r;\lambda)}{r}=0, \mbox{ for any } 
0< \lambda \leq 1. \]
The growth of 
$f_X(l,r;\lambda)$ is not a quasi-isometry invariant 
since the modulus of the annulus might be changed by a quasi-isometry.  
However the property of having a weak sub-linear filling area 
is a quasi-isometry invariant, as it can be easily checked. 
It is not difficult to see that the proof of our theorem actually 
shows that a manifold of bounded geometry and 
linear volume growth whose filling growth is weakly sublinear 
should be simple connected at infinity.

 \bibliographystyle{plain}

\end{document}